\documentclass[twoside]{article} 
\date{} 
\title{Asymptotics of some integrals involving modified Bessel and hyper-Bessel functions}
\author{\sc R. B.\ Paris \\
{\em Division of Computing and Mathematics,} \\
{\em Abertay University, Dundee DD1 1HG, UK}}
\begin{document}
\def\f#1#2{\mbox{${\textstyle \frac{#1}{#2}}$}}
\def\dfrac#1#2{\displaystyle{\frac{#1}{#2}}}
\def\boldal{\mbox{\boldmath $\alpha$}}
\newcommand{\bee}{\begin{equation}}
\newcommand{\ee}{\end{equation}}
\newcommand{\lam}{\lambda}
\newcommand{\ka}{\kappa}
\newcommand{\al}{\alpha}
\newcommand{\la}{\lambda}
\newcommand{\ga}{\gamma}
\newcommand{\sa}{\sigma}
\newcommand{\eps}{\epsilon}
\newcommand{\fr}{\frac{1}{2}}
\newcommand{\fs}{\f{1}{2}}
\newcommand{\g}{\Gamma}
\newcommand{\br}{\biggr}
\newcommand{\bl}{\biggl}
\newcommand{\ra}{\rightarrow}
\newcommand{\gtwid}{\raisebox{-.8ex}{\mbox{$\stackrel{\textstyle >}{\sim}$}}}
\newcommand{\ltwid}{\raisebox{-.8ex}{\mbox{$\stackrel{\textstyle <}{\sim}$}}}
\renewcommand{\topfraction}{0.9}
\renewcommand{\bottomfraction}{0.9}
\renewcommand{\textfraction}{0.05}
\newcommand{\mcol}{\multicolumn}
\date{}
\maketitle
\pagestyle{myheadings}
\markboth{\hfill \sc R. B.\ Paris  \hfill}
{\hfill \sc Asymptotic expansion of an integral\hfill}
\begin{abstract}
We investigate the asymptotic expansion of integrals analogous to Ball's integral
\[\int_0^\infty \bl(\frac{\g(1+\nu)|J_\nu(x)|}{(x/2)^\nu}\br)^{\!n}dx\]
for large $n$ in which the Bessel function $J_\nu(x)$ is replaced by the modified Bessel functions $I_\nu(x)$ and $K_\nu(x)$ together with appropriate exponential factors $e^{\mp x}$, respectively.

The above integral with $J_\nu(x)$ replaced by a hyper-Bessel function of the type recently discussed in
Aktas {\it et al.} [The Ramanujan J., 2019] and taken over a finite interval determined by the first positive zero of the function is also considered for $n\to\infty$. We give the leading asymptotic behaviour of the hyper-Bessel function for $x\to+\infty$ in an appendix. Numerical examples are given to illustrate the accuracy of the various expansions obtained.
\vspace{0.3cm}

\vspace{0.1cm}
 
\end{abstract}

\vspace{0.3cm}

\noindent $\,$\hrulefill $\,$

\vspace{0.3cm}

\begin{center}
{\bf 1.\ Introduction}
\end{center}
\setcounter{section}{1}
\setcounter{equation}{0}
\renewcommand{\theequation}{\arabic{section}.\arabic{equation}}
The asymptotic expansion of Ball's integral \cite{KB} for large positive values of $n$
\bee\label{e10}
\int_0^\infty \bl(\frac{\g(1+\nu) |J_\nu(x)|}{(x/2)^\nu}\br)^{\!n} x^{2\nu-1}dx,\qquad n\geq 2,\ \nu\geq \fs,
\ee
where $J_\nu(x)$ is the Bessel function of the first kind, was investigated by Kerman {\it et al.} in \cite{KOS}. In a recent note \cite{P} it was shown that the above integral could be replaced by
\bee\label{e100}
\int_0^{j_{\nu,1}} \bl(\frac{\g(1+\nu) J_\nu(x)}{(x/2)^\nu}\br)^{\!n} x^{2\nu-1}dx
\ee
to within exponentially small terms when $n$ is large, where $j_{\nu,1}$ is the first positive zero of $J_\nu(x)$. The large-$n$ expansion was found in the form
\[2^{2\nu-1} (1+\nu)^\nu \g(\nu)\sum_{k=0}^\infty \frac{(-)^k c_k}{n^{k+\nu}}\qquad (n\to\infty),\]
where the leading coefficient $c_0=1$; explicit values of $c_k$ for $k\leq 3$ were obtained in \cite{KOS} and for $k\leq 6$ in \cite{P}. 

Following a suggestion of T. Pog\'any \cite{Pog}, we consider analogous integrals involving the modified Bessel functions $I_\nu(x)$ and $K_\nu(x)$ and also the hyper-Bessel function $J_{\sa_1, \ldots ,\sa_m}(x)$ defined in Section 4. 
In Section 2, we consider the expansion of the integral
\[{\cal I}_n=\int_0^\infty \bl(e^{-x} \frac{\g(1+\nu) I_\nu(x)}{(x/2)^\nu}\br)^{\!n} dx\qquad(\nu>-\fs)\]
for $n\to\infty$. Since $I_\nu(x)\sim e^x/\sqrt{2\pi x}$ as $x\to+\infty$, it is necessary to add the factor $e^{-x}$ to cancel the exponential growth of $I_\nu(x)$. The integrand is then of O($x^{-\nu-1/2})$ as $x\to\infty$ so that ${\cal I}_n$ converges for $\nu>-\fs$. In Section 3, a similar process is adopted to determine the large-$n$ expansion of
\[{\cal K}_n=\int_0^\infty \bl(\frac{2e^x}{\g(\nu)} (x/2)^\nu K_\nu(x)\br)^{\!-n}dx\qquad (\nu>\fs),\]
where $n>0$.
From the small and large argument behaviours
\[K_\nu(x)\sim \frac{1}{2}\g(\nu) (x/2)^{-\nu} \quad (x\to 0), \qquad K_\nu(x)\sim \sqrt{\frac{\pi}{2x}}\,e^{-x}\quad (x\to+\infty),\]
it is seen that the integrand has the value unity at $x=0$ and is O($x^{\nu-1/2})$ as $x\to\infty$, thereby necessitating the condition $\nu>\fs$ for convergence.

In the final section, we consider an integral analogous to (\ref{e10}) in which the classical Bessel function $J_\nu(x)$ is replaced by a hyper-Bessel function. A significant difference, however, is that the interval of integration cannot be taken as $[0,\infty)$ on account of the asymptotic structure of the particular hyper-Bessel function under consideration. It is necessary to take a finite integration interval analogous to that in (\ref{e100}) determined by the first positive zero of the function.

\vspace{0.6cm}

\begin{center}
{\bf 2.\ An integral involving the modified Bessel function $I_\nu(x)$}
\end{center}
\setcounter{section}{2}
\setcounter{equation}{0}
\renewcommand{\theequation}{\arabic{section}.\arabic{equation}}
The first integral we consider is the analogue of (\ref{e10}) where the function $|J_\nu(x)|$ is replaced by the modified Bessel function $I_\nu(x)$, viz.
\bee\label{e31}
{\cal I}_n=\int_0^\infty \bl(e^{-x} \frac{\g(1+\nu) I_\nu(x)}{(x/2)^\nu}\br)^{\!n} dx\qquad (\nu>-\fs),
\ee
where $n>0$ (not necessarily an integer). The integrand has the value unity at $x=0$ and is a monotonically decreasing function. This can be seen by
letting $y(x)=e^{-x} (\fs x)^{-\nu} I_\nu(x)$, whence
\[y'(x)=-e^{-x} (\fs x)^{-\nu} \{I_\nu(x)-I_{\nu+1}(x)\}<0,\]
since for $x>0$ and $\nu>-\fs$ it is known  that $I_\nu(x)>I_{\nu+1}(x)$ (see \cite{J}).

We write
\[\psi(x)=x-\log\,\bl(\frac{\g(1+\nu) I_\nu(x)}{(x/2)^\nu}\br)=x-\log\,\sum_{k=0}^\infty \frac{(x/2)^{2k}}{(1+\nu)_k k!}.\]
Since $\psi(0)=0$ and $\psi(\infty)=\infty$, the change of variable $\tau=\psi(x)$ yields 
\bee\label{e32}
{\cal I}_n=\int_0^\infty e^{-n\tau}\frac{dx}{d\tau}\,d\tau,
\ee
where 
\[\tau=\psi(x)=x-\frac{x^2}{4(1+\nu)}+\frac{x^4}{32(1+\nu)^2(2+\nu)}-\frac{x^6}{96(1+\nu)^3(2+\nu)(3+\nu)}+ \cdots\]
valid\footnote{The circle of convergence is determined by the nearest singularity of $\psi(x)$ that occurs at $x=\pm ij_{\nu,1}$, since $I_\nu(z)$ on the imaginary $z$-axis behaves like $J_\nu(|z|)$.} in $x<j_{\nu,1}$.
Inversion of this series with the help of {\it Mathematica} then produces
\[x=\tau+\frac{\tau^2}{4(1+\nu)}+\frac{\tau^3}{8(1+\nu)^2}+\frac{(8+3\nu)\tau^4}{64(1+\nu)^3 (2+\nu)}+\frac{(8+\nu)\tau^4}{128(1+\nu)^4(2+\nu)}+\cdots,\]
whence we obtain the expansion
\bee\label{e33}
\frac{dx}{d\tau}=\sum_{k=0}^\infty A_k \tau^k \qquad (\tau<\tau_0).
\ee
The first few coefficients $A_k$ are
\begin{eqnarray*}
A_0&=&1,\quad A_1=\frac{1}{2(1+\nu)},\quad A_2=\frac{3}{8(1+\nu)^2},\\
A_3&=&\frac{8+3\nu}{16(1+\nu)^3 (2+\nu)},\quad A_4=\frac{5(8+\nu)}{128(1+\nu)^4 (2+\nu)},\quad A_5=\frac{142+11\nu-5\nu^2}{256(1+\nu)^5(2+\nu)(3+\nu)},\\ A_6&=&\frac{7(272+8\nu-71\nu^2-5\nu^3)}{3072(1+\nu)^6(2+\nu)^2(3+\nu)},\quad
A_7=\frac{2656-1364\nu-1602\nu^2-137\nu^3+19\nu^4}{2048(1+\nu)^7(2+\nu)^2(3+\nu)(4+\nu)},\\
A_8&=&\frac{3(6816-11740\nu-4770\nu^2+1025\nu^3+119\nu^4)}{32768(1+\nu)^8(2+\nu)^2(3+\nu)(4+\nu)}, \ldots \ .
\end{eqnarray*}
The expansion (\ref{e33}) holds in $\tau<\tau_0$, where $\tau_0=|\psi(\pm ij_{\nu,2}')|$ since $x=\pm ij_{\nu,2}'$ is the nearest point in the mapping $x\mapsto \tau$ where $dx/d\tau$ is singular. The quantity $j_{\nu,2}'$ is the second positive zero of $J_\nu'(x)$.

From (\ref{e32}) and (\ref{e33}), straightforward integration yields
\[{\cal I}_n\sim \sum_{k=0}^\infty A_k \int_0^\infty e^{-n\tau} \tau^k d\tau=\sum_{k=0}^\infty \frac{k! A_k}{n^{k+1}}.\]
Then we obtain:
\newtheorem{theorem}{Theorem}
\begin{theorem}$\!\!\!.$\ \ For $\nu>-\fs$ and $n\to\infty$ the following expansion holds
\[
{\cal I}_n\sim \frac{1}{n}\bl\{1+\frac{1}{2(1+\nu)n}+\frac{3}{4(1+\nu)^2 n^2}+ \frac{3(8+3\nu)}{8(1+\nu)^3(2+\nu)n^3}+
\frac{15(8+\nu)}{16(1+\nu)^4(2+\nu)n^4}\]
\bee\label{e33th}
+\frac{15(142+11\nu-5\nu^2)}{32(1+\nu)^5(2+\nu)(3+\nu)n^5}+\frac{105(272+8\nu-71\nu^2-5\nu^3)}{64(1+\nu)^6(2+\nu)^2(3+\nu) n^6}+ \cdots \br\}.
\ee
\end{theorem}
\bigskip

\begin{table}[t]
\caption{\footnotesize{Values of the absolute relative error in the computation of ${\cal I}_n$  when $n=100$ for different values of the order $\nu$ and truncation index $k$ using the expansion in (\ref{e33th}).}}
\begin{center}
\begin{tabular}{|r|r|r|r|}
\hline
&&&\\[-0.35cm]
\mcol{1}{|c|}{$k$} & \mcol{1}{c|}{$\nu=0$} & \mcol{1}{c|}{$\nu=3/4$} & \mcol{1}{c|}{$\nu=1$} \\
[.1cm]\hline
&&&\\[-0.35cm]
0 & $5.051\times 10^{-03}$ & $2.874\times 10^{-03}$ & $2.513\times 10^{-03}$\\
1 & $7.615\times 10^{-05}$ & $2.468\times 10^{-05}$ & $1.888\times 10^{-05}$\\
2 & $1.531\times 10^{-06}$ & $2.633\times 10^{-07}$ & $1.732\times 10^{-07}$\\
3 & $3.845\times 10^{-08}$ & $3.213\times 10^{-09}$ & $1.772\times 10^{-09}$\\
4 & $1.142\times 10^{-09}$ & $4.119\times 10^{-11}$ & $1.817\times 10^{-11}$\\
5 & $3.841\times 10^{-11}$ & $4.762\times 10^{-13}$ & $1.444\times 10^{-13}$\\
6 & $1.409\times 10^{-12}$ & $2.327\times 10^{-15}$ & $5.266\times 10^{-16}$\\
[.1cm]\hline
\end{tabular}
\end{center}
\end{table}
\noindent In Table 1 we show values of the absolute relative error in the computation of ${\cal I}_n$ against truncation index
for different values of $\nu$.

A similar integral is given by
\[{\hat{\cal I}}_n=\int_0^\infty \bl(\frac{\g(1+\nu) I_\nu(x)}{(x/2)^\nu}\br)^{\!-n}dx \qquad (\nu\geq 0),\]
where $n>0$. With the standard substitution $\psi(x)=\log\,(\g(1+\nu)I_\nu(x)/(x/2)^\nu)$ and change of variable $\tau^2=\psi(x)$, we find
\[\tau^2=\frac{x^2}{4(1+\nu)}-\frac{x^4}{32(1+\nu)^2(2+\nu)}+\frac{x^6}{96(1+\nu)^3(2+\nu)(3+\nu)}+\cdots\ ,\]
which upon inversion yields the expansion
\[\frac{1}{2\sqrt{1+\nu}} \,\frac{dx}{d\tau}=\sum_{k=0}^\infty {\hat A}_k \tau^{2k}\qquad (\tau<\tau_0).\]
The first few coefficients ${\hat A}_k$ are
\begin{eqnarray*}
{\hat A}_0\!\!&=&\!\!1,\quad {\hat A}_1=\frac{3}{4(2+\nu)},\quad {\hat A}_2=-\frac{5(1+11\nu)}{96(2+\nu)^2(3+\nu)},\\
{\hat A}_3\!\!&=&\!\!-\frac{7(20+9\nu-17\nu^2)}{128(2+\nu)^3(3+\nu)(4+\nu)},\\
{\hat A}_4\!\!&=&\!\!\frac{75404+262439\nu+182205\nu^2-28031\nu^3-19409\nu^4}{10240(2+\nu)^4(3+\nu)^2(4+\nu)(5+\nu)},\\
{\hat A}_5\!\!&=&\!\!\frac{11(127864-364742\nu-1417421\nu^2-966731\nu^3+8605\nu^4+48361\nu^5}{122880(2+\nu)^5(3+\nu)^2(4+\nu)(5+\nu)(6+\nu)}.
\end{eqnarray*}
We note that in this case the integrand has a saddle point at $x=0$. Routine evaluation then produces the asymptotic expansion
\bee\label{e34}
{\hat{\cal I}}_n\sim (1+\nu)^{1/2} \sum_{k=0}^\infty \frac{{\hat A}_k \g(k+\fs)}{n^{k+1/2}}\qquad (n\to\infty).
\ee

\vspace{0.6cm}

\begin{center}
{\bf 3.\ An integral involving the modified Bessel function $K_\nu(x)$}
\end{center}
\setcounter{section}{3}
\setcounter{equation}{0}
\renewcommand{\theequation}{\arabic{section}.\arabic{equation}}
In this section we consider an analogous integral to (\ref{e31}) with the modified Bessel function of the second kind, also known as the Macdonald function $K_\nu(x)$, namely
\bee\label{e41}
{\cal K}_n=\int_0^\infty \bl(\frac{2e^x}{\g(\nu)} (x/2)^\nu K_\nu(x)\br)^{\!-n}dx\qquad (\nu>\fs),
\ee 
where again $n>0$ (not necessarily an integer) and we assume throughout this section that $\nu$ is not an integer.  The quantity in brackets in (\ref{e41}) is monotonically increasing for $\nu>\fs$, since with $y(x)=x^\nu e^x K_\nu(x)$ we have
\[y'(x)=x^\nu e^x\{K_\nu(x)-K_{\nu-1}(x)\}>0.\]
The fact that the quantity in braces is positive follows from the result \cite[(10.32.9)]{DLMF}
\[K_\nu(x)-K_{\nu-1}(x)=\int_0^\infty e^{-x\cosh t} \{\cosh \nu t-\cosh (\nu-1)t\}\,dt>0,\qquad \nu>\fs.\]

With the changes of variable $\psi(x)=x+\log\,\{2(x/2)^\nu K_\nu(x)/\g(\nu)\}$ and $\tau=\psi(x)$, so that $x\in[0,\infty)$ maps to $\tau\in[0,\infty)$ when $\nu>\fs$, we have
\[{\cal K}_n=\int_0^\infty e^{-n\psi(x)} dt=\int_0^\infty e^{-n\tau}\,\frac{dx}{d\tau}\,d\tau.\]
To proceed further we require the inversion of the mapping $\tau\mapsto x$. To do this in general terms is complicated so we prefer to carry out this procedure for specific values of $\nu$. 
From the definition
\[K_\nu(x)=\frac{\pi}{2\sin \pi\nu} \{I_{-\nu}(x)-I_\nu(x)\},\]
it is seen that, provided $\nu\neq 1, 2, \ldots\ $,
\[\frac{2}{\g(\nu)} (x/2)^\nu K_\nu(x)=1+\frac{x^2}{4(1-\nu)}+\frac{x^4}{32(1-\nu)(2-\nu)}+\cdots\hspace{4cm}\] \[\hspace{4cm}+\bl(\frac{x}{2}\br)^{\!2\nu}\frac{\g(-\nu)}{\g(\nu)}\bl(1+\frac{x^2}{4(1+\nu)}+\frac{x^4}{32(1+\nu)(2+\nu)}+\cdots\br).\]
Thus when $\nu=2/3$, for example, we find
\[\tau=x+g \bl(\frac{x}{2}\br)^{\!4/3}+\frac{3}{4}x^2-g^2 \bl(\frac{x}{2}\br)^{\!8/3}-\frac{12}{5}g \bl(\frac{x}{2}\br)^{\!10/3}+\bl(-\frac{27}{128}+\frac{1}{48}g^3\br)x^4+\cdots\ ,\]
where $g:=\g(-2/3)/\g(2/3)$, which upon inversion yields
\[x=\tau-\frac{g\tau^{4/3}}{2^{4/3}}+\frac{g^2 \tau^{5/3}}{3\cdot 2^{2/3}}+(\frac{g^3}{8}-\frac{3}{4})\,\tau^2+\frac{5g}{648\cdot 2^{1/3}} (13g^3+162) \tau^{7/3}+\cdots\ \]
valid in a neighbourhood of $\tau=0$. 

However, it is found much easier to deal with this inversion process with the coefficients expressed in numerical form rather than in algebraic form. In this manner we obtain after differentiation with respect to $\tau$
\[\frac{dx}{d\tau}=\sum_{k=0}^\infty B_k(\f{2}{3}) \tau^{k/3},\]
where the coefficients $B_k(\f{2}{3})$ are listed in Table 2 for $k\leq 6$. Then
\bee\label{e42}
{\cal K}_n\sim \sum_{k=0}^\infty \frac{B_k(\f{2}{3})}{n^{k/3+1}} \int_0^\infty e^{-w} w^{k/3}dw=
\sum_{k=0}^\infty \frac{B_k(\f{2}{3})}{n^{k/3+1}}\,\g(k/3+1)\qquad (\nu=\f{2}{3})
\ee
as $n\to\infty$.

We present the series expansion for $dx/d\tau$ for the two cases  $\nu=6/5$ and  $\nu=4/3$. The coefficients $B_k(\nu)$ are computed using the two lines of {\it Mathematica} commands below and are given in Table 2. 
\begin{eqnarray*}
&&\tt{f:=2(x/2)^\nu BesselK[\nu,x]/Gamma(\nu);S= N[Series[x+Log\,[f],\{x,0,m\}],18]}\\
&&\tt{D[InverseSeries[S,\tau],\tau]}
\end{eqnarray*}
where $m$ is an integer that determines how far we carry out the expansion process.
The asymptotic expansion of ${\cal K}_n$ is then computed as above.
\begin{table}[h]
\caption{\footnotesize{The coefficients $B_k(\nu)$ for $1\leq k\leq 6$ (with $B_0(\nu)=1$) for different values of the order $\nu$.}}
\begin{center}
\begin{tabular}{|r|r|r|r|}
\hline
&&&\\[-0.35cm]
\mcol{1}{|c|}{$k$} & \mcol{1}{c|}{$\nu=2/3$} & \mcol{1}{c|}{$\nu=6/5$} & \mcol{1}{c|}{$\nu=4/3$} \\
[.1cm]\hline
&&&\\[-0.35cm]
1 & $+ 1.570228753470$ & $+2.5000000000$ & $+1.5000000000$\\
2 & $+ 3.082022922779$ & $-2.4023937306$ & $-1.4329122397$\\
3 & $+ 5.033299366471$ & $+9.3750000000$ & $+3.3750000000$\\
4 & $+ 7.536861569983$ & $-18.7186511511$& $-6.8958901536$\\
5 & $+ 10.658049516385$& $+9.1382014252$ & $+10.1250000000$\\
6 & $+ 14.372332247788$& $+42.9687500000$& $+3.3365109159$\\
[.1cm]\hline
\end{tabular}
\end{center}
\end{table}
When $\nu=6/5$, we have
\[\frac{dx}{d\tau}=1+B_1(\f{6}{5}) \tau+B_2(\f{6}{5})\tau^{7/5}+B_3(\f{6}{5})\tau^2+B_4(\f{6}{5})\tau^{12/5}+B_5(\f{6}{5})\tau^{14/5}+B_6(\f{6}{5})\tau^3+\cdots\]
and the asymptotic expansion
\bee\label{e43}
{\cal K}_n\sim\frac{1}{n}\bl\{1+\frac{B_1(\f{6}{5})}{n}+\frac{B_2(\f{6}{5})}{n^{7/5}}\g(\f{12}{5})+\frac{2B_3(\f{6}{5})}{n^2}+\frac{B_4(\f{6}{5})}{n^{12/5}} \g(\f{17}{5})+\frac{B_5(\f{6}{5})}{n^{16/5}} \g(\f{19}{5})+\frac{6B_6(\f{6}{5})}{n^3}+\cdots\br\};
\ee
when $\nu=4/3$, we have
\[\frac{dx}{d\tau}=1+B_1(\f{4}{3}) \tau+B_2(\f{4}{3})\tau^{5/3}+B_3(\f{4}{3})\tau^2+B_4(\f{4}{3})\tau^{8/3}+B_5(\f{4}{3})\tau^{3}+B_6(\f{4}{3})\tau^{10/3}+\cdots\]
and the asymptotic expansion
\bee\label{e44}
{\cal K}_n\sim\frac{1}{n}\bl\{1+\frac{B_1(\f{4}{3})}{n}+\frac{B_2(\f{4}{3})}{n^{5/3}}\g(\f{8}{3})+\frac{2B_3(\f{4}{3})}{n^2}+\frac{B_4(\f{4}{3})}{n^{8/3}} \g(\f{11}{3})+\frac{6B_5(\f{4}{3})}{n^{3}}+\frac{B_6(\f{4}{3})}{n^{10/3}} \g(\f{13}{3})+\cdots\br\}
\ee
as $n\to\infty$.

In Table 3 we present the absolute relative error in the computation of the integral ${\cal K}_n$  for the three values of the order $\nu$ using different truncations of the expansions in (\ref{e42})--(\ref{e44}). Because these expansions involve inverse fractional powers of $n$, it is seen that the rate of decay of the relative error with increasing truncation index is rather slow.
\begin{table}[t]
\caption{\footnotesize{Values of the absolute relative error in the computation of ${\cal K}_n$  when $n=100$ for different values of the order $\nu$ and truncation index $k$ using the expansions in (\ref{e42})--(\ref{e44}).}}
\begin{center}
\begin{tabular}{|r|r|r|r|}
\hline
&&&\\[-0.35cm]
\mcol{1}{|c|}{$k$} & \mcol{1}{c|}{$\nu=2/3$} & \mcol{1}{c|}{$\nu=6/5$} & \mcol{1}{c|}{$\nu=4/3$} \\
[.1cm]\hline
&&&\\[-0.35cm]
0 & $3.391\times 10^{-1}$ & $2.104\times 10^{-2}$ & $1.439\times 10^{-2}$\\
1 & $1.394\times 10^{-1}$ & $3.435\times 10^{-3}$ & $3.931\times 10^{-4}$\\
2 & $5.407\times 10^{-2}$ & $1.195\times 10^{-3}$ & $5.932\times 10^{-4}$\\
3 & $2.080\times 10^{-2}$ & $6.402\times 10^{-4}$ & $7.201\times 10^{-5}$\\
4 & $8.025\times 10^{-3}$ & $2.256\times 10^{-4}$ & $5.450\times 10^{-5}$\\
5 & $3.106\times 10^{-3}$ & $1.201\times 10^{-4}$ & $5.371\times 10^{-6}$\\
[.1cm]\hline
\end{tabular}
\end{center}
\end{table}

\vspace{0.6cm}

\begin{center}
{\bf 4.\ An integral involving the hyper-Bessel function}
\end{center}
\setcounter{section}{4}
\setcounter{equation}{0}
\renewcommand{\theequation}{\arabic{section}.\arabic{equation}}
The particular hyper-Bessel function we shall use to replace the Bessel function $J_\nu(x)$ in the integral (\ref{e10}) is defined by \cite{ABS}
\bee\label{e11}
J_{\sa_1, \ldots ,\sa_m}(x)=\frac{(x/(m+1))^{\sa_1+\cdots +\sa_m}}{\prod_{j=1}^m \g(\sa_j+1)}\,{}_0F_m\bl(\!\!\begin{array}{c}-\!\!-\\ \sa_1\!+\!1, \ldots , \sa_m\!+\!1\end{array}\!;-\bl(\frac{x}{m+1}\br)^{\!m+1}\br).
\ee
Here ${}_0F_m$ denotes the generalised hypergeometric function with $m$ denominator parameters
\[{}_0F_m(z)=\sum_{k=0}^\infty \frac{1}{(\sa_1\!+\!1)_k \ldots (\sa_m\!+\!1)_k}\,\frac{z^k}{k!}\qquad (\sa_j>-1,\ \ 1\leq j\leq m)
\]
and $(a)_k=\g(a+k)/\g(a)=a(a+1)\ldots (a+k-1)$ is Pochhammer's symbol for the rising factorial.
When $m=1$, $\sa_1=\nu$, the definition (\ref{e11}) reduces to the classical Bessel function $J_\nu(x)$, viz.
\[J_\nu(x)=\frac{(\fs x)^\nu}{\g(1+\nu)}\,{}_0F_1\bl(\!\!\begin{array}{c}-\!\!-\\ 1+\nu\end{array}\!;-\bl(\frac{ x}{2}\br)^{\!2}\br)=\frac{(\fs x)^\nu}{\g(1+\nu)}\,\sum_{k=0}^\infty \frac{(-)^k (\fs x)^{2k}}{(1+\nu)_k k!}~.\]

Before we can formulate an integral analogous to that in (\ref{e10}), it is necessary to consider the basic properties and asymptotic behaviour of $J_{\sa_1, \ldots ,\sa_m}(x)$. 
In what follows we write ${\vec \sigma}=\{\sa_1, \sa_2, \ldots , \sa_m\}$ and define the quantity
\bee\label{e12}
\mu_k:=\prod_{j=1}^m (\sa_j+k)^{-1}.
\ee
Both $J_{\sa_1,\ldots , \sa_m}(x)$ (for $m\geq 2$) and $J_\nu(x)$ have an infinite number of zeros on $[0, \infty)$.
If the first such zero of $J_{\sa_1, \ldots ,\sa_m}(x)$ is denoted by $j_{\vec\sa,1}$, it was established in \cite[Theorem 4]{ABS} that
\[(m+1) \mu_1^{-1/(m+1)} <j_{\vec\sa,1} < (m+1) (\mu_1-\mu_2)^{-1/(m+1)}.\]
However, although these two functions possess similar zero properties, their asymptotic structure is quite different.
From (\ref{a1}) in the appendix, the leading asymptotic behaviour of $J_{\sa_1, \ldots ,\sa_m}(x)$ is 
\[J_{\sa_1, \ldots ,\sa_m}(x)\sim\frac{2(2\pi)^{-m/2}}{(m+1)^{1/2}} \bl(\frac{x}{m+1}\br)^{\!-m/2}\,e^{x \cos \pi/(m+1)}\,\cos \br\{x \sin \frac{\pi}{m+1}+\frac{\pi\vartheta}{m+1}\br\}\]
as $x\to+\infty$, where $\vartheta$ is defined in (\ref{a0}). 

Thus, when $m\geq 2$ the hyper-Bessel function grows exponentially as $x\to+\infty$ 
and consequently the integral (\ref{e10}) modified to incorporate the hyper-Bessel function cannot be taken over an infinite range. Accordingly, we consider the asymptotic expansion of the integral over the finite interval $[0, j_{\vec\sa,1}]$ viz.
\bee\label{e13}
{\cal J}_n=\int_0^{j_{\vec\sa,1}} \bl({\cal J}_{\sa_1, \ldots , \sa_m}(x)\br)^{\!n}dx,
\ee
where the normalised hyper-Bessel function ${\cal J}_{\sa_1,\ldots ,\sa_m}(x)$  is defined by
\begin{eqnarray}
{\cal J}_{\sa_1,\ldots ,\sa_m}(x)&=&\frac{\prod_{j=1}^m\g(\sa_j+1)}{(x/(m+1))^{\sa_1+\cdots+\sa_m}}\,J_{\sa_1,\ldots ,\sa_m}(x)\nonumber\\
&=&{}_0F_m\bl(\!\!\begin{array}{c}-\!\!-\\ \sa_1\!+\!1, \ldots , \sa_m\!+\!1\end{array}\!;-\bl(\frac{x}{m+1}\br)^{\!m+1}\br).\label{e15}
\end{eqnarray}
In the case $m=1$, $\sa_1=\nu$, the integral (\ref{e13}) reduces to that in (\ref{e100}) when the factor $x^{2\nu-1}$ is omitted. An equivalent factor could be added to (\ref{e13}), but we choose not to do this in order to avoid the appearance of additional parameters.

\vspace{0.3cm}

\noindent{\bf 4.1\ \ The asymptotic expansion}
\vspace{0.3cm}

\noindent
Let $p=m+1$ and set $\psi(x)=-\log\,{}_0F_m(-(x/p)^p)$, so that the integral (\ref{e13}) becomes
\[{\cal J}_n=\int_0^{j_{\vec\sa,1}} e^{-n\psi(x)}dx.\]
Since $\psi(0)=0$ and $\psi(j_{\vec\sa,1})=\infty$, then with the change of variable $\tau^p=\psi(x)$ we have
\bee\label{e21}
{\cal J}_n=\int_0^\infty e^{-n\tau^p} \,\frac{dx}{d\tau}\,d\tau,
\ee
where
\begin{eqnarray*}
\tau^p&=&\psi(x)=-\log\,\bl(1-\frac{(x/p)^p \mu_1}{1!}+\frac{(x/p)^{2p} \mu_1 \mu_2}{2!}-\frac{(x/p)^{3p} \mu_1 \mu_2 \mu_3}{3!}+\cdots \br)\\
&=&\mu_1 (x/p)^p+\frac{1}{2}\mu_1(\mu_1-\mu_2) (x/p)^{2p}+\frac{1}{6}(2\mu_1^3-3\mu_1^2 \mu_2+\mu_1 \mu_2 \mu_3) (x/p)^{3p}+\cdots\ 
\end{eqnarray*}
valid in $x<p (j_{\vec\sa,1})^{1/p}$. Inversion of this last expression with the aid of {\it Mathematica} yields
\[(x/p)^p=\frac{\tau^p}{\mu_1}+\frac{(\mu_2-\mu_1)}{2\mu_1^2}\,\tau^{2p}+ \frac{(\mu_1^2-3\mu_1 \mu_2+3\mu_2^2-\mu_2\mu_3)}{6\mu_1^3}\,\tau^{3p}+\cdots  \]
whence
\[x=\frac{p\tau}{\mu_1^{1/p}}\bl\{1+\frac{(\mu_2-\mu_1)}{2\mu_1}\,\tau^{p}+\frac{(\mu_1^2-3\mu_1 \mu_2+3\mu_2^2-\mu_2\mu_3)}{6\mu_1^2}\,\tau^{2p} +\cdots \br\}^{\!1/p}.\]
This then leads to an expansion for $dx/d\tau$ given by
\bee\label{e22}
\frac{dx}{d\tau}=\frac{p}{\mu_1^{1/p}} \sum_{k=0}^\infty \frac{(-)^k A_k }{p^k}\, (kp+1)\tau^{kp}\qquad (\tau<\tau_0),
\ee
where
\begin{eqnarray}
A_0\!\!\!\!&=&\!\!\!\!1,\quad A_1=\frac{1}{2}(1-\ga_2),\nonumber\\
A_2\!\!\!\!&=&\!\!\!\!\frac{1}{24}\bl\{p+3-6(p+1)\ga_2+3(3p+1)\ga_2^2-4p\ga_2 \ga_3\br\},\nonumber\\
A_3\!\!\!\!&=&\!\!\!\!\frac{1}{48}\bl\{p+1-(p+1)(4p+3)\ga_2+3(2p+1)(3p+1)\ga_2^2-4p(2p+1)\ga_2\ga_3\br\},\nonumber\\
A_4\!\!\!\!&=&\!\!\!\!\frac{1}{5760}\bl\{15\! +\! 30 p\! +\! 5 p^2\! -\! 2 p^3\!-\! 60 (p\! +\! 1)^2 (2 p\!+\!1)\ga_2 + 
 10 (3\! +\! 13 p\! +\! 14 p^2) \ga_2 (3(3p\! +\! 1) \ga_2 - 4 p \ga_3) \nonumber\\
 &&\!\!\!\!\!\!- 
 60 (3p\! +\! 1) \ga_2 ((1 \!+\! 9 p\! +\! 20 p^2) \ga_2^2 - 4 p (4p\!+\!1) \ga_2\ga_3 + 
    2 p^2 \ga_3 \ga_4) + \ga_2 \{15 (1\! +\! 18 p\! +\! 107 p^2\! +\! 210 p^3) \ga_2^3\nonumber\\
    &&\!\!\!\!\!\! - 
    120 p (1\! +\! 11 p\! +\! 30 p^2) \ga_2^2 \ga_3 + 
    40 p^2 (5 p\!+\!1) \ga_2 \ga_3 (2 \ga_3\! +\! 3 \ga_4) - 48 p^3 \ga_3\ga_4\ga_5\}\br\}\label{e22c}
\end{eqnarray}
with
\[\ga_k:=\frac{\mu_k}{\mu_1}=\prod_{j=1}^m \frac{\sa_j+1}{\sa_j+k}\qquad (k\geq 2).\]
The expansion (\ref{e22}) holds in $\tau<\tau_0$, where $\tau_0^p=|\psi(j_{\vec\sa,2}')|$ since $x=j_{\vec\sa,2}'$ is the nearest point in the mapping $x\mapsto\tau$ where $dx/d\tau$ is singular. The quantity $j_{\vec\sa,2}'$ is the second positive zero of the derivative of the hyper-Bessel function, which interlaces with the zeros $j_{\vec\sa,1}$ and $j_{\vec\sa,2}$ \cite[Theorem 5]{ABS}. 

Then from (\ref{e21}) and (\ref{e22}), we obtain
\begin{eqnarray*}
{\cal J}_n&\sim&\frac{p}{\mu_1^{1/p}} \sum_{k=0}^\infty \frac{(-)^k A_k}{p^k}\,(kp+1) \int_0^\infty e^{-n\tau^p} \tau^{kp}\,d\tau\\
&=&\frac{1}{(n\mu_1)^{1/p}} \sum_{k=0}^\infty \frac{(-)^k A_k}{(np)^k}\,(kp+1) \int_0^\infty e^{-w} w^{k+1/p-1} dw.
\end{eqnarray*}
Evaluation of the integral as a gamma function then produces
\begin{theorem}$\!\!.$\ \ With $p=m+1$ and $\mu_1=\prod_{j=1}^m (\sa_j+1)^{-1}$, we have the expansion
\bee\label{e23}
{\cal J}_n\sim \frac{p}{(n\mu_1)^{1/p}} \sum_{k=0}^\infty \frac{(-)^k A_k}{(np)^k} \,\g(k+\frac{1}{p}+1)
\ee
as $n\to\infty$, where the first five coefficients $A_k$ are listed in (\ref{e22c}). 
\end{theorem}

In Table 4 we show values of the absolute relative error against truncation index $k$ in the evaluation of ${\cal J}_n$ in the case $m=2$ using (\ref{e23}). The value of ${\cal J}_n$ was obtained by a high-precision numerical integration procedure
with the first zero when $\sa_1=\fs$, $\sa_2=\f{3}{4}$ given by $j_{\vec\sa,1}\doteq 4.5079255667$.
\begin{table}[th]
\caption{\footnotesize{Values of the absolute relative error in the computation of ${\cal J}_n$ against truncation index $k$ when $m=2$ and $\sa_1=1/2$, $\sa_2=3/4$. .}}
\begin{center}
\begin{tabular}{|r|l|l|l|}
\hline
&&&\\[-0.35cm]
\mcol{1}{|c|}{$k$} & \mcol{1}{c|}{$n=20$} & \mcol{1}{c|}{$n=50$} & \mcol{1}{c|}{$n=100$}\\
[.1cm]\hline
&&&\\[-0.35cm]
0 & $6.916\times 10^{-03}$ & $2.753\times 10^{-03}$ & $1.375\times10^{-03}$ \\
1 & $1.260\times 10^{-05}$ & $1.819\times 10^{-06}$ & $4.387\times10^{-07}$ \\
2 & $1.953\times 10^{-06}$ & $1.233\times 10^{-07}$ & $1.534\times10^{-08}$ \\
3 & $3.243\times 10^{-08}$ & $8.988\times 10^{-10}$ & $5.753\times10^{-11}$ \\
4 & $4.606\times 10^{-09}$ & $4.532\times 10^{-11}$ & $1.396\times10^{-12}$ \\
[.1cm]\hline
\end{tabular}
\end{center}
\end{table}
\vspace{0.6cm}

\begin{center}
{\bf Appendix: \ The asymptotic behaviour of $J_{\sa_1,\ldots , \sa_m}(x)$ for $x\to+\infty$}
\end{center}
\setcounter{section}{1}
\setcounter{equation}{0}
\renewcommand{\theequation}{\Alph{section}.\arabic{equation}}
The hypergeometric-type function
\[f(z)=\sum_{k=0}^\infty \frac{z^k}{\prod_{j=1}^m \g(\sa_j+1+k)k!}\]
is associated with the parameters\footnote{Here we follow the notation of \cite[\S 2.3]{PK} and denote by $\kappa$ the quantity $m+1$, although in Section 4 this quantity was denoted by $p$} (see \cite[\S 2.3]{PK})
\bee\label{a0}
\kappa=m+1,\quad h=1,\quad \vartheta=-\frac{1}{2}m+\sum_{j=1}^m \sa_j.
\ee
Define the formal exponential asymptotic sum 
\[E(z):=Z^\vartheta e^Z \sum_{k=0}^\infty A_k Z^{-k}, \qquad Z:=\kappa(hz)^{1/\kappa},\]
where $A_k$ are constants independent of $z$ with $A_0=(2\pi)^{-m/2} \kappa^{-\fr-\vartheta}$.
Then,  when $\kappa>2$ (that is, $m\geq 2$) the asymptotic expansion of $f(z)$ is given by \cite[\S 2.3, Case (iii)]{PK}
\[f(z)\sim \sum_{r=-P}^P E(z e^{2\pi ir})\qquad (|z|\to\infty,\ |\arg\,z|\leq\pi),\]
where $P$ is chosen such that $2P+1$ is the smallest odd integer satisfying $2P+1>\fs\kappa$.

For the hypergeometric function appearing in (\ref{e11}) we have $Z=xe^{\pi i/\kappa}$. Then, when $m\geq 2$,
\[{}_0F_m(-(x/\kappa)^\kappa)\sim, \prod_{j=1}^m\g(\sa_j+1)\,\sum_{r=-P}^P E(xe^{(2r+1)\pi i})\qquad (x\to+\infty).\]
The dominant exponential sums correspond to $r=0$ and $r=-1$, whence we obtain
\begin{eqnarray*}
{}_oF_m(-(x/\kappa)^\kappa)&\sim& \prod_{j=1}^m\g(\sa_j+1)\,\{E(xe^{\pi i})+E(xe^{-\pi i})\}\\
&\sim& 2A_0\prod_{j=1}^m\g(\sa_j+1)\,x^\vartheta e^{x\cos\,\pi/\kappa} \cos \bl\{x \sin\,\frac{\pi}{\kappa}+\frac{\pi\vartheta}{\kappa}\br\}.
\end{eqnarray*}
Hence the laeding behaviour of $J_{\sa_1, \ldots , \sa_m}(x)$ is given by
\bee\label{a1}
J_{\sa_1, \ldots ,\sa_m}(x)\sim\frac{2(2\pi)^{-m/2}}{(m+1)^{1/2}} \bl(\frac{x}{m+1}\br)^{\!-m/2}\,e^{x \cos \pi/(m+1)}\,\cos \br\{x \sin \frac{\pi}{m+1}+\frac{\pi\vartheta}{m+1}\br\}
\ee
as $x\to+\infty$; see also \cite{PW}.

When $m=1$, $\sa_1=\nu$, the approximation (\ref{a1}) reduces to the well-known leading behaviour of the classical Bessel function \cite[(10.17.3)]{DLMF}
\[J_\nu(x)\sim \sqrt{\frac{2}{\pi x}}\,\cos\bl\{x-\frac{\pi\nu}{2}-\frac{\pi}{4}\br\} \qquad (x\to+\infty).\]
However, when $m\geq 2$ it is seen from (\ref{a1}) that $J_{\sa_1, \ldots , \sa_m}(x)$ is oscillatory with an exponentially growing amplitude as $x\to+\infty$, and so is of a completely different asymptotic structure to that of $J_\nu(x)$.

\end{document}